\documentclass{svmult}
\usepackage{graphicx} 

\usepackage{epsfig}
\usepackage{amssymb}
\usepackage{bm}
\usepackage{epstopdf}
\usepackage{amsmath}

\begin{document}
\title*{\begin{center}Robust Estimation of Location \\ through Schoenberg transformations \end{center}}
\toctitle{Robust Estimation through Schoenberg transformations}
%
\titlerunning{Robust Estimation through Schoenberg transformations}
%
\author{
  Fran\c{c}ois Bavaud\inst{}
}

\institute{
Department of Computer Science and Mathematical Methods\\
Department of Geography\\
University of Lausanne,  Switzerland\\
\date{September 2010}}

\maketitle             

\begin{abstract}
Schoenberg transformations, mapping  Euclidean configurations into Euclidean configurations, 
 define in turn a transformed inertia, whose minimization produces robust location estimates. The procedure only depends upon Euclidean distances between observations, and applies equivalently to univariate and multivariate data. The choice of the family of transformations and their parameters defines a  flexible location strategy,   generalizing $M$-estimators. 
Two regimes of solutions are identified. Theoretical results on their  existence and stability are provided, and illustrated on two data sets.
\end{abstract}

 {\bf Keywords:} Correspondence Analysis, Euclidean distances, Huber function, Huygens principles, $M$-estimators, Schoenberg transformations, Tukey bisquare

\section{Introduction and main result}
This paper investigates the properties of presumably new location estimates, defined through 
Schoenberg transformations, which map initial Euclidean distances into new Euclidean distances (Schoenberg 1938; Bavaud 2011).

Specifically, consider an univariate sample $\{x_i\}_{i=1}^n$  with weights  $\{f_i\}_{i=1}^n$. The weighted mean  $a=\bar{x}_f:=\sum_i f_i x_i$ minimizes $\sum_i f_i (x_i-a)^2$, and the weighted median $a=x_{0.5}$ minimizes $\sum_i f_i |x_i-a|$. 
This paper  introduces and studies a  family of {\em centroids} $a$, suitable for any dimension, generalizing  the previous well-known basic estimates. The centroids are defined as the averages
$a=\sum_i \alpha_i x_i$ minimizing the {\em transformed inertia} 
\begin{equation}
\label{funct}
\Gamma(a)=\sum_{i=1}^n f_i \varphi(D_{ia})
\end{equation}
where $D_{ia}=\|x_i-a\|^2$ is a squared Euclidean distance and $\varphi(D)$ a {\em Schoenberg transformation}   (see section \ref{scchch}).

The profile $\alpha$ generating minimizers $a=\sum_i \alpha_i x_i$ turns out to satisfy
\begin{equation}
\label{solmain}
\alpha_i =\frac{f_i \varphi'(D_{ia})}{\sum_j f_j \varphi'(D_{ja})}
\qquad\qquad\qquad
D_{ia}=\sum_j \alpha_j D_{ij}-\frac12\sum_{jk}\alpha_j\alpha_k D_{jk}\enspace .
\end{equation}
This pair of identities defines an iterative scheme, depending on the distances between observations only, and converging in general towards a local minimum of   (\ref{funct}). The term $\varphi'(D_{ia})$ downweights observations distant from the centroid $a$, which behaves as a multidimensional robust estimate of location,  
comparable  to the $M$-estimates in the one-dimensional case.

Section  \ref{deap}  defines the main ingredients and presents general results, in particular the existence of two regimes of solutions. Section  \ref{univ} shows the connection with the theory of $M$-estimates in one dimension, illustrated by the  
 {\tt copper}  {\tt concentration} data in Section \ref{copper}. Section \ref{multi} illustrates 
 the multidimensional case by means of the  chi-square distances occurring in correspondence analysis.

\section{Definitions and general results}
\label{deap}
Data are characterized by the matrix of  dissimilarities $D_{ij}$ between observations $i,j=1,\ldots,n$, together with their weights $f_i>0$ with $\sum_i f_i=1$. We assume the dissimilarities to be squared Euclidean, that is of the form $D_{ij}=\|x_i-x_j\|^2$ where the coordinates $x_i\in \mathbb{R}^p$, unique up to a translation and a rotation,  can be recovered by MDS (see e.g. Mardia et al. 1979; Borg and Groenen 1997; 
Bavaud 2011; and references therein), in a space of dimension $p\le n-1$. 

Here and in the sequel, {\em we assume that observations are distinct}, that is $D_{ij}>0$ for $i\neq j$. That is, possibly identical initial observations $x_{i}=x_j$ should be first aggregated into a single observation of weight $f_i+f_j$.

\subsection{Huygens principles}
\label{huyhuy}
The {\em inertia} $\Delta(f)$,  measuring the dispersion of the weighted configuration, 
expresses  as in either equivalent two forms  (Huygens weak principle)
\begin{equation}
\label{inetertrte}
\Delta(f)=\frac12\sum_{i,j=1}^nf_i f_j D_{ij}=\sum_i f_i D_{if}
\end{equation}
where 
$D_{if}=\|x_i-\bar{x}_f\|^2$ and 
$\bar{x}_f =\sum_i f_i x_i$. Also (Huygens strong principle)
\begin{equation}
\label{inetertrte2}
\sum_j f_j D_{ij}=D_{if}+\Delta(f)\enspace .
\end{equation}
Consider another distribution or profile $\alpha$ with $\sum_i\alpha_i=1$, with associated centroid $a=
\sum_i\alpha_i x_i$. Substituting $\alpha$ to $f$ in (\ref{inetertrte2}) yields the second identity in 
(\ref{solmain}). Also, $\sum_i f_i D_{ia}=D_{fa}+\Delta(f)$, which shows $\Gamma(a)$ in 
(\ref{funct}) to be minimum for $a=\bar{x}_f $ in the identity case  $\varphi(D)=D$ - an elementary result.

\subsection{Schoenberg transformations}
\label{scchch}
A  {\em Schoenberg transformation} is  a  componentwise mapping $\tilde{D}_{ij}= \phi(D_{ij})$ with the property that if $D_{ij}$ represents a squared Euclidean distance between observations $i$ and $j$, so does  $\tilde{D}_{ij}$ (irrespectively of the dimension $p$). 

The class of all Schoenberg transformations has been investigated and determined by Schoenberg (1938, Theorem 6 p. 828; Bavaud 2011):

\begin{theorem}[Schoenberg 1938]
\label{schoen}
The function $\phi(D)$ is a Schoenberg transformation iff   of the form 
\begin{equation}
\label{schtr}
\varphi(D)=\int_0^\infty \frac{1-\exp(-\lambda D)}{\lambda}\: g(\lambda)\: d\lambda
\end{equation}
where $g(\lambda) \:  d\lambda$ is a non-negative measure on $[0,\infty)$ such that $\int_{1}^\infty \frac{g(\lambda)}{\lambda}  d\lambda<\infty$. 

\

Equivalently, $\varphi(D)$ is a Schoenberg transformation iff it is smooth, with $\varphi(0)=0$, 
$\varphi^{(2r-1)}(D)\ge0$ and $\varphi^{(2r)}(D)\le0$ for all $r=1,2,\ldots$
\end{theorem}
The second part is a consequence of Bernstein theorem on completely monotonic functions. In particular, a Schoenberg transformation is increasing, concave,  and zero at the origin. Examples are provided by (Bavaud 2011):
\begin{eqnarray}
\label{poexlo}
\begin{array}{lllcr}
\phi(D)&=&D^q\quad&(0<q<1) \quad & \mbox{power transformation}\\
\phi(D)&=& 1-\exp(-  D/\delta))\quad&(\delta>0)  \quad & \mbox{exponential transformation}\\
\phi(D)&=& \ln(1+D/\delta)\quad&(\delta>0)   \quad & \mbox{logarithmic transformation.}\\
 \end{array}\end{eqnarray}

A transformation $\varphi(D)$ is said to be  {\it rectifiable} if $\varphi'(0)<\infty$, that is iff  $\int_0^\infty  g(\lambda)  \: d\lambda<\infty$. By (\ref{schtr}), rectifiable transformations obtain as mixtures of $(1-\exp(-\lambda D)/\lambda$, which tends to the identity  $\tilde{D}=D$ for  $\lambda\to0$.

A transformation $\varphi(D)$  is said to be
{\it bounded} if $\varphi(\infty)<\infty$, that is iff  $\int_0^\infty   \frac{g(\lambda)}{\lambda} \: d\lambda<\infty$. By (\ref{schtr}), bounded transformations obtain as mixtures of $1-\exp(-\lambda D)$, which tends for  $\lambda\to\infty$ to the discrete metric $\tilde{D}=I(D>0)$,  attributing a unit 
dissimilarity between distinct observations. 

 The power transformation is not rectifiable nor bounded; the exponential transform is rectifiable and bounded; the logarithmic transformation is rectifiable but not bounded. 
 
\subsection{High-dimensional embedding and strain}
 \label{hdes}
To the $n\times n$ matrix of transformed squared distances $\tilde{D}_{ij}$ correspond, by MDS, new coordinates $\tilde{x}_i\in \mathbb{R}^{\tilde{p}}$, unique up to a rotation and  translation. Thus, 
a  Schoenberg transformation induces a mapping or {\em embedding}
$\tilde{x}=\eta(x)$ of the 
original  coordinates into the transformed  coordinates, similar  to  the  high-dimensional embeddings of Machine Learning (Bavaud 2011). 

\

Let $\tilde{a}$ denote a position in the the transformed space, and define 
\begin{displaymath}
\tau(\tilde{a})=\sum_i f_i  \|\tilde{x}_i-\tilde{a}\|^2\enspace. 
\end{displaymath}
On one hand, 
\begin{displaymath}\min_{\tilde{a}} \tau(\tilde{a})=\frac12\sum_{ij}f_if_j\tilde{D}_{ij}=
\tilde{\Delta}(f)
\end{displaymath}
in view of Section \ref{huyhuy}. On the other hand, minimizing the transformed inertia 
$\Gamma(a)=\tau(\eta(a))$ amounts in minimizing $ \tau(\tilde{a})$ {\em under the additional constraint 
$\tilde{a}=\eta(a)=\eta(\sum_i \alpha_i x_i)$} for some $\alpha$.  The importance of that constraint, reflecting  the non-linearity of the embedding $\eta$, can be measured 
 by the {\em strain} $\gamma(a)\triangleq   \Gamma(a)/\tilde{\Delta}(f)$, obeying $\gamma(a)
\ge1$ by construction.

\subsection{Behavior of minima}
Two types of minima exist:  {\em distributed minima}, where the identities  (\ref{solmain}) hold, and $\alpha$ is strictly positive at each observation;  and {\em concentrated minima}, where $\alpha$ is concentrated on some observation $i_0$, thus making $a=x_{i_0}$. The latter case holds iff the quantity  $\max_i  \phi'(D_{ia})$ is infinite.

\begin{theorem}
\label{tworegimes}

$\mbox{}$ 

1a)  when the transformation is rectifiable, a minimizer  $a=\sum_i\alpha_i x_i$ of $\Gamma(a)$  necessarily satisfies the identities in  (\ref{solmain}) (distributed case)

1b) the minimizer also necessarily satisfies the stability condition
\begin{equation}
\label{scdo}
\min_b\sum_i f_i [\varphi'(D_{ia})  +2\varphi''(D_{ia})D_{ia}\cos^2\theta_{ib}^a]\: \ge\: 0
\end{equation} 
where $b$ is another point and $\theta_{ib}^a$ the angle between the vectors $x_i-a$  and $b-a$. In particular, $a$ is stable if
\begin{equation}
\label{stabily2bis}
\sum_i f _i [\phi'(D_{ia})+2\phi''(D_{ia})D_{ia}]\ge0\enspace .
\end{equation}

2) when the transformation is not rectifiable, a minimizer either behaves  as in the   distributed case 1), or  in a concentrated way, with support concentrated on a single  observation $x_{i_0}$.
\end{theorem}

\begin{proof}
See the appendix for  a proof of 1a) and 1b). To prove 2), note that equation (\ref{solmain}) is not justified anymore iff $\varphi'(D_{ia})=\infty$ for some $i$. As $\varphi'(D)$ is always finite for $D>0$, then, necessarily, $D_{ia}=0$ and $\varphi'(0)=\infty$ (non-rectifiable
 transformation). The observations being distinct, this situation depicts a centroid concentrated on an unique $i=i_0$ (concentrated case). In particular, this regime necessarily holds whenever 
the distributed stability condition (\ref{scdo}) is violated. $\Box$
\end{proof}

The   {\em concentration}  of the minimizing profile  $\alpha$ can be measured by its {\em entropy} $H(\alpha)\triangleq -\sum_j \alpha_j \ln \alpha_j\ge0$, where $H(\alpha)=0$ iff the profile is  concentrated.

\section{The univariate case}
\label{univ}
In one dimension,   minimizing   $\Gamma(a)=
\sum_i f_i\:  \phi((x_i-a)^2)$ yields
\begin{displaymath}
0=\sum_i f_i\:  \phi'((x_i-a)^2)(x_i-a)\equiv\sum_i f_i \: \psi(x_i-a)\enspace . 
\end{displaymath}
The odd function defined for $D\ge0$ as 
\begin{equation}
\label{psifunc}
\psi(\sqrt{D})\triangleq\phi'(D)\sqrt{D}
\end{equation}
is known as the ``$\psi$-function" in the  theory of  $M$-estimators (e.g. 
Huber 1964, Hampel et al. 1986, Maronna et al. 2006 and references therein). It is fair to add that the iterative scheme (\ref{solmain}) also generalizes the  $W$-estimation of location proposed by Tukey (1977),  itself identified as a  particular  case of the $M$-estimation (Hampel et al. 1986 p.116).

 The sign of the derivative 
\begin{equation}
\label{chidef}
\chi(D)\triangleq \frac{\partial \psi(\sqrt{D})}{\partial\sqrt{D}}=\phi'(D)+2\phi''(D)D
\end{equation}
governs the stability of the $\phi$-estimate, in accordance with Theorem \ref{tworegimes} 1b), where
conditions  (\ref{scdo})  and  (\ref{stabily2bis}) now coincide in view of $\cos\theta_{ib}^a=\pm1$ in one dimension. 

\

\subsection{Power transformation}
\label{potrans}
The {\em power transformation} $\phi(D)=D^q$ with $0< q< 1$ is not rectifiable nor bounded. It   exhibits both regimes, namely, distributed solutions for $q>1/2$ and concentrated ones for $q<1/2$, as demonstrated by the study of the sign of $\chi(D)=q(2q-1)D^{q-1}$ (see figure \ref{figpower}).

\subsection{Rectifiable transformations}
Transformations satisfying $\varphi'(0)<\infty$ can be written as  
\begin{equation}
\label{retrcl}
\varphi(D)=\eta(\delta)\: h(\frac{D}{\delta})
\qquad\qquad h(0)=0\qquad\qquad h'(0)=1
\end{equation}
where $\delta$ is a (squared) characteristic length, and $\eta(\delta)$ is immaterial for our purpose. Such are  the exponential and logarithmic transformations of Section \ref{scchch}, as well as the {\em Tukey} and {\em Huber} transformations
\begin{eqnarray}
\label{tuhu}
\begin{array}{lllr}
\varphi(D)&=& 
\left\{\begin{array}{lcr} D-D^2/\delta+D^3/(3\delta^2)     &\quad \mbox{if }   D\le \delta \\ 1/3  &   \quad\mbox{otherwise} \end{array}\right.  & \hspace{0.2cm}\mbox{ ``Tukey" transformation}\\
 \varphi(D)&=&
\left\{\begin{array}{lcr} D     & \quad\hspace{1.4cm} \mbox{ if }   D\le \delta \\ 2\sqrt{\delta D}-\delta      &   \quad\hspace{1.9cm}\mbox{otherwise} \end{array}\right.  &\hspace{0.3cm}\mbox{ ``Huber" transformation.}\\
\end{array}\end{eqnarray}
As a matter of fact, the $\psi$-functions  associated to those transformations  can be shown to respectively yield the so-called  {\em Tukey bisquare} and the {\em Huber function}, familiar in the theory of $M$-estimates, whence their names. The derivatives  of the Tukey  and the Huber transformations  follow the pattern of Theorem \ref{schoen},  except at the value $D=\delta$, presenting  a discontinuity in the third, respectively second-order derivative.

\begin{theorem}[rectifiable transformations; proof in the appendix]
\label{trecti}
There exists a finite characteristic length $\delta_0$ granting the uniqueness of the minimizer $a$ for $\delta\ge\delta_0$, with $\lim_{\delta \to \infty}a=
 \bar{x}_f$. \end{theorem}

\subsection{Bounded transformations}
Bounded transformations can be written as    
\begin{equation}
\label{boncl}
\varphi(D)=\eta(\delta)\: h(\frac{D}{\delta})
\qquad\qquad h(0)=0\qquad\qquad h(\infty)=1\enspace .
\end{equation}
Such is the case of the discrete metric $\varphi(D)=I(D>0)$,  attributing a unit dissimilarity between distinct observations. The transformed inertia reads $\Gamma(a)=1-\sum_i f_i I(a=x_i)$, and  attains its minimum values $1-f_i$  when $a=x_i$, thus generating concentrated solutions only.

This behavior is emblematic of bounded transformations in general (rectifiable or not), which tend  towards  the discrete metric in the limit $\delta\to0$ (Section \ref{scchch}). Hence, by continuity, $n$ minima emerge in the 
 limit of small characteristic length:

\begin{theorem}[bounded transformations]
\label{trsddsti}
There exists a finite characteristic length $\delta_0$ granting the existence of $n$ minima, 
  for $\delta\le\delta_0$,  each being  located near an observation $x_i$.
  \end{theorem}

\begin{figure}
\begin{center}
\includegraphics[width=5.8cm]{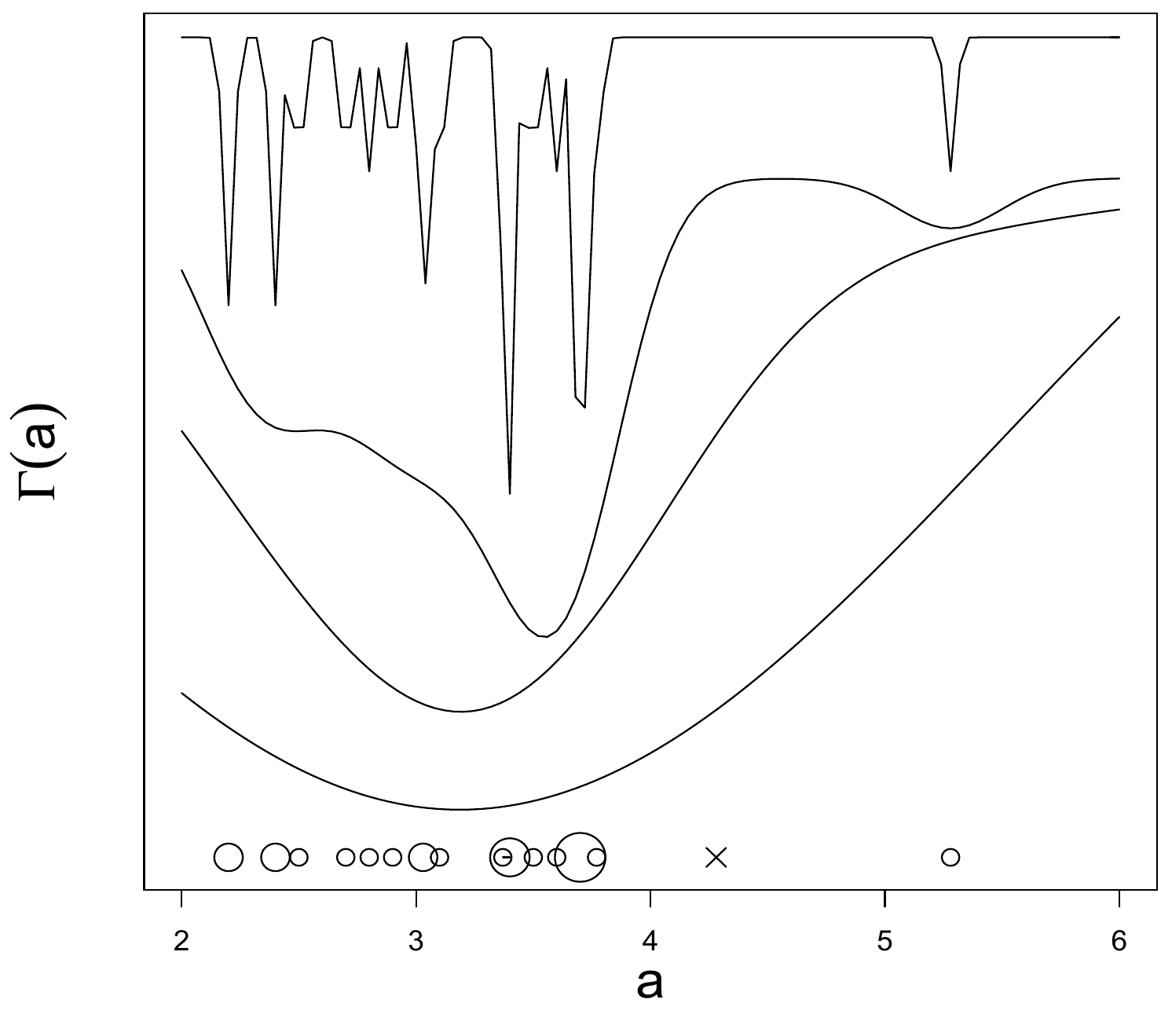}
\includegraphics[width=5.8cm]{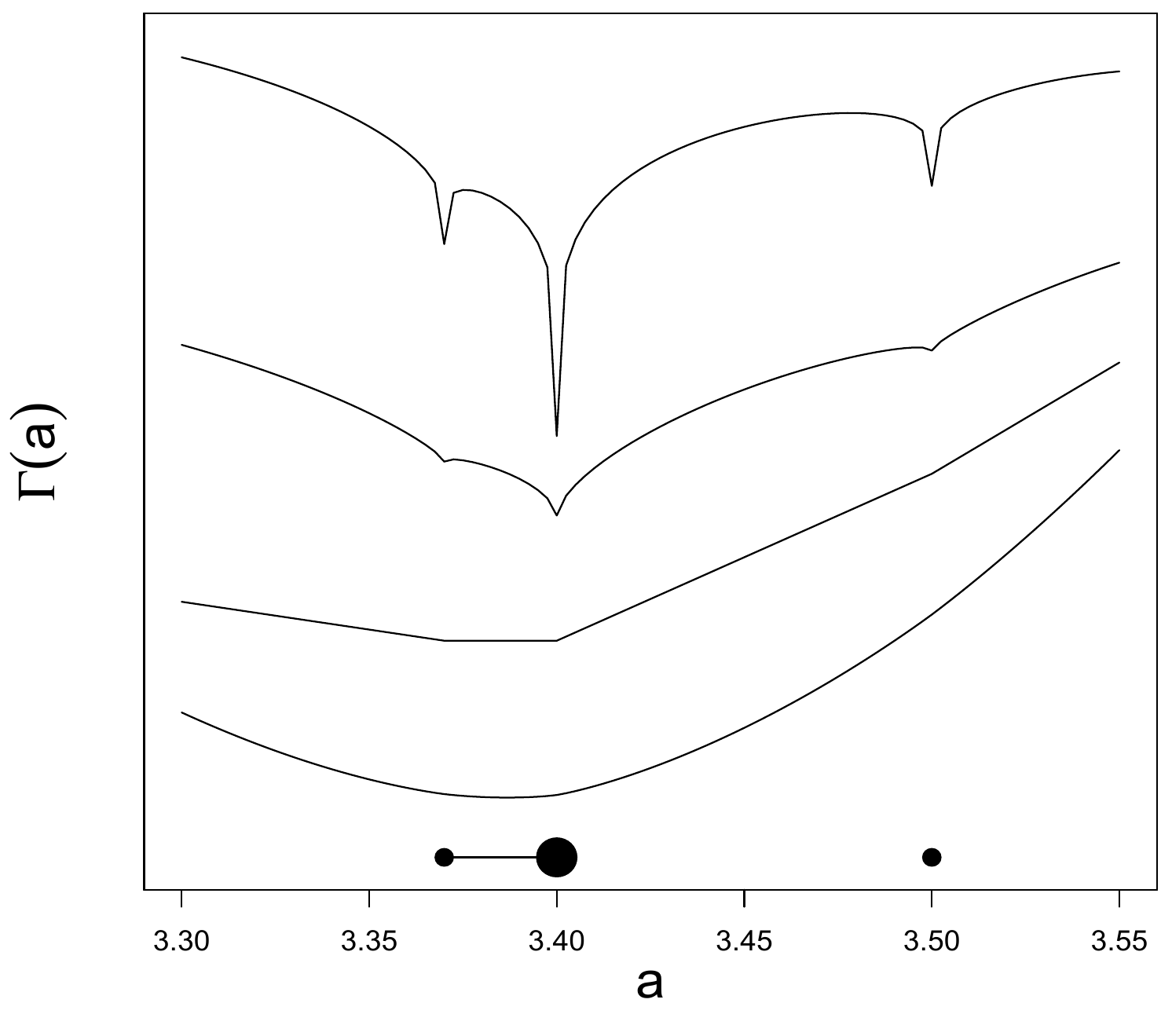}
\caption{Transformed inertia $\Gamma(a)$ (arbitrary units) versus $a$. Left: exponential transform $\varphi(D)=1-\exp(-\lambda D)$ with   $\lambda$ =0.1, 1, 10 and 1'000, in ascending order. Right: power transform $\varphi(D)=D^q$ with $q$ =0.7, 0.5, 0.3 and 0.1, in ascending order. The bottom line depicts the  values of 15, respectively 3 weighted observations. ``$\times$" denotes $\bar{x}_f$, and the segment denotes the median interval.}
\label{objectiveplot}
\end{center}
\end{figure}

 \section{Illustration:  copper concentrations data}
 \label{copper}
Consider the dataset {\tt chem}, available in the R library {\tt MASS} (Venables and Ripley  2002) consisting of $n=24$ copper concentrations in wholemeal flour (Abbey, 1988), with sorted values 

\

\noindent {\tt 2.20 2.20 2.40 2.40 2.50 2.70 2.80 2.90 3.03
3.03 3.10 3.37}
 
\noindent  {\tt   3.40 3.40 3.40 3.50 3.60 3.70
3.70 3.70 3.70 3.77 5.28 28.95}Ê

\

Ties occur 
twice ($x$ = 2.20, 2.40, 3.03), three times ($x$ = 3.40) or four times ($x$ = 3.70), and must
preliminarily aggregated,   resulting in a sample of $n=16$ observations with varying  weights.

\

Figure \ref{objectiveplot} plots the functions $\Gamma(a)$ resulting form the exponential and the power transform, for varying parameters. It illustrates the gradual emergence of $n$ minima from the mean $\bar{x}_f=4.2804$, in accordance to Section \ref{univ}. In the power case at the transition value $q=1/2$, the  minimum is attained at the {\em median interval} $x_{0.5}=(3.37,3.40)$ (second plot of figure \ref{objectiveplot}, right). 
 
\

\begin{figure}
\begin{center}
\includegraphics[width=2.1in]{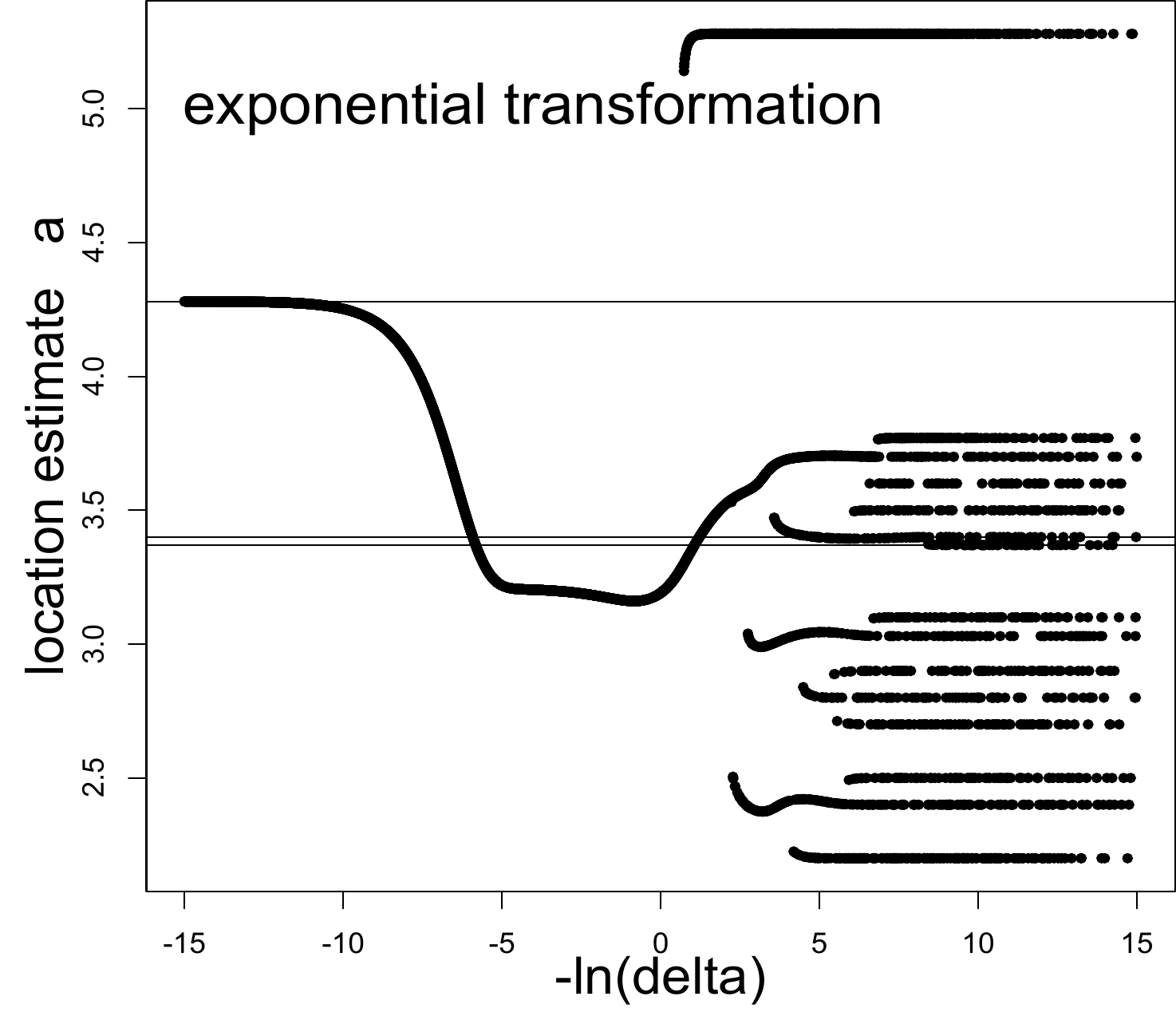}\hspace{0.5cm}
\includegraphics[width=2.1in]{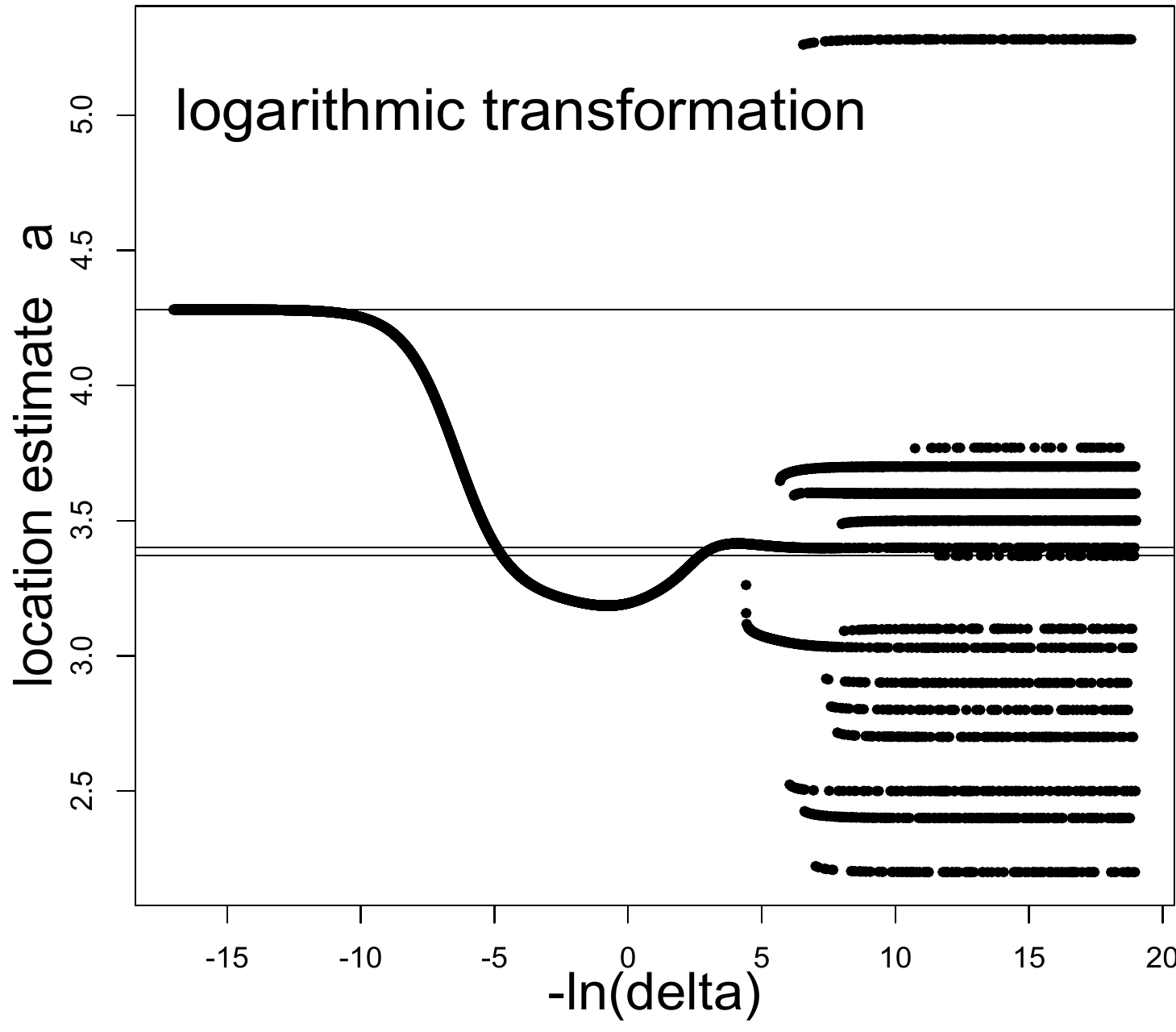}
\includegraphics[width=2.1in]{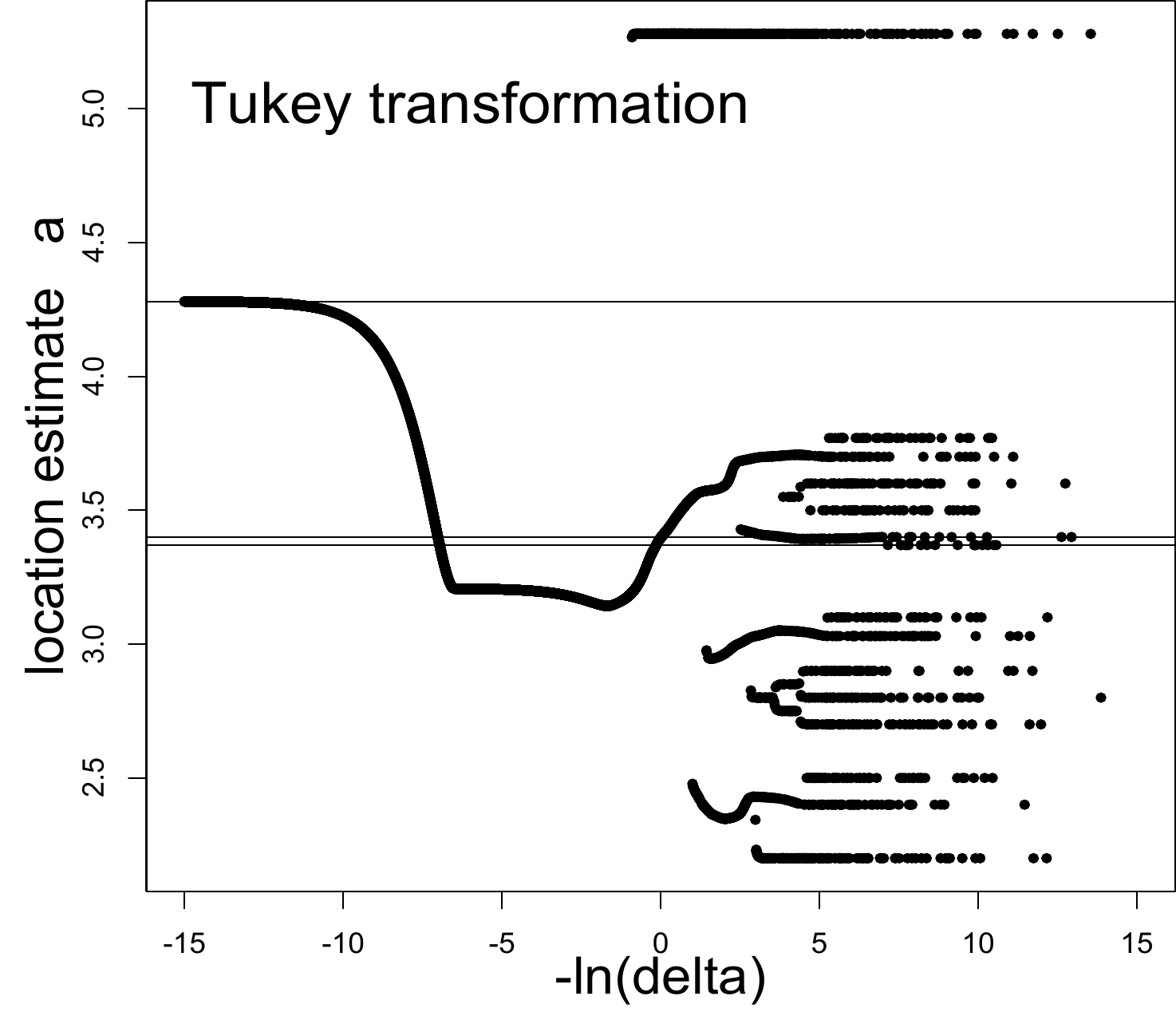}\hspace{0.5cm}
\includegraphics[width=2.1in]{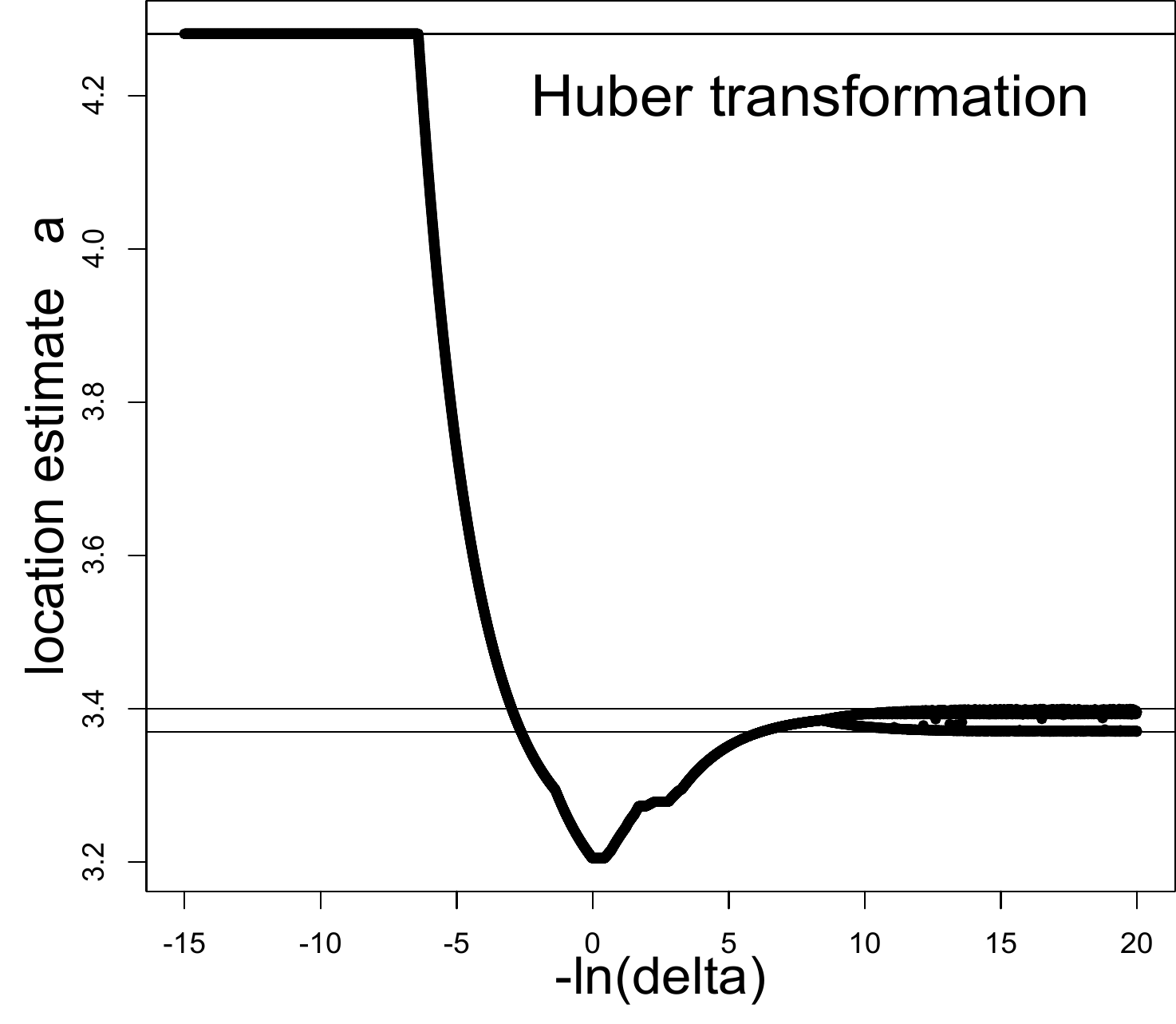}
\caption{Location estimates for the copper concentrations data, under various rectifiable transforms: exponential and logarithmic (\ref{poexlo}), Tukey and Huber   (\ref{tuhu}).}
\label{fourtransforms}
\end{center}
\end{figure}

\begin{figure}
\begin{center}
\includegraphics[width=2.1in]{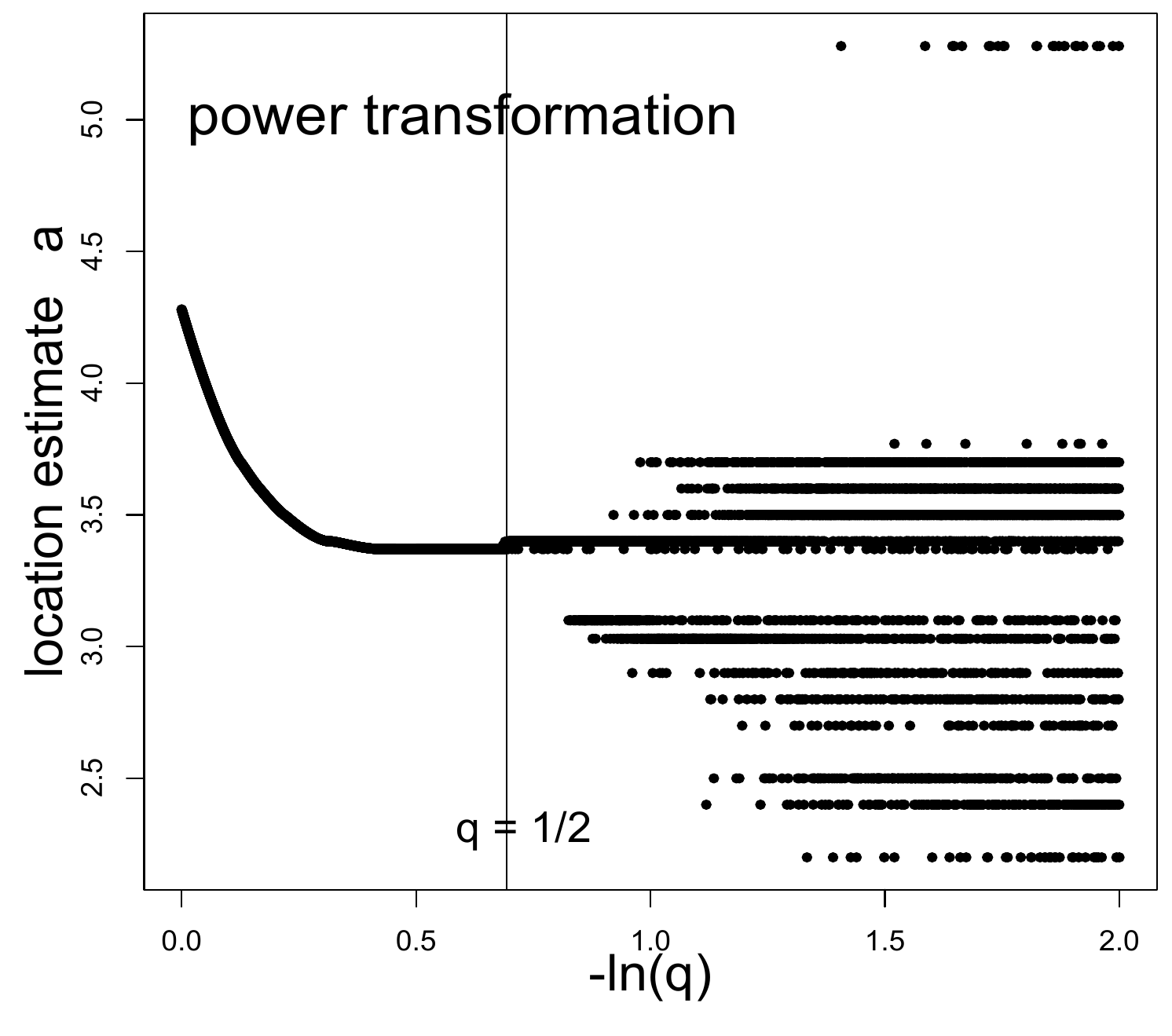}\hspace{0.5cm}
\includegraphics[width=2.1in]{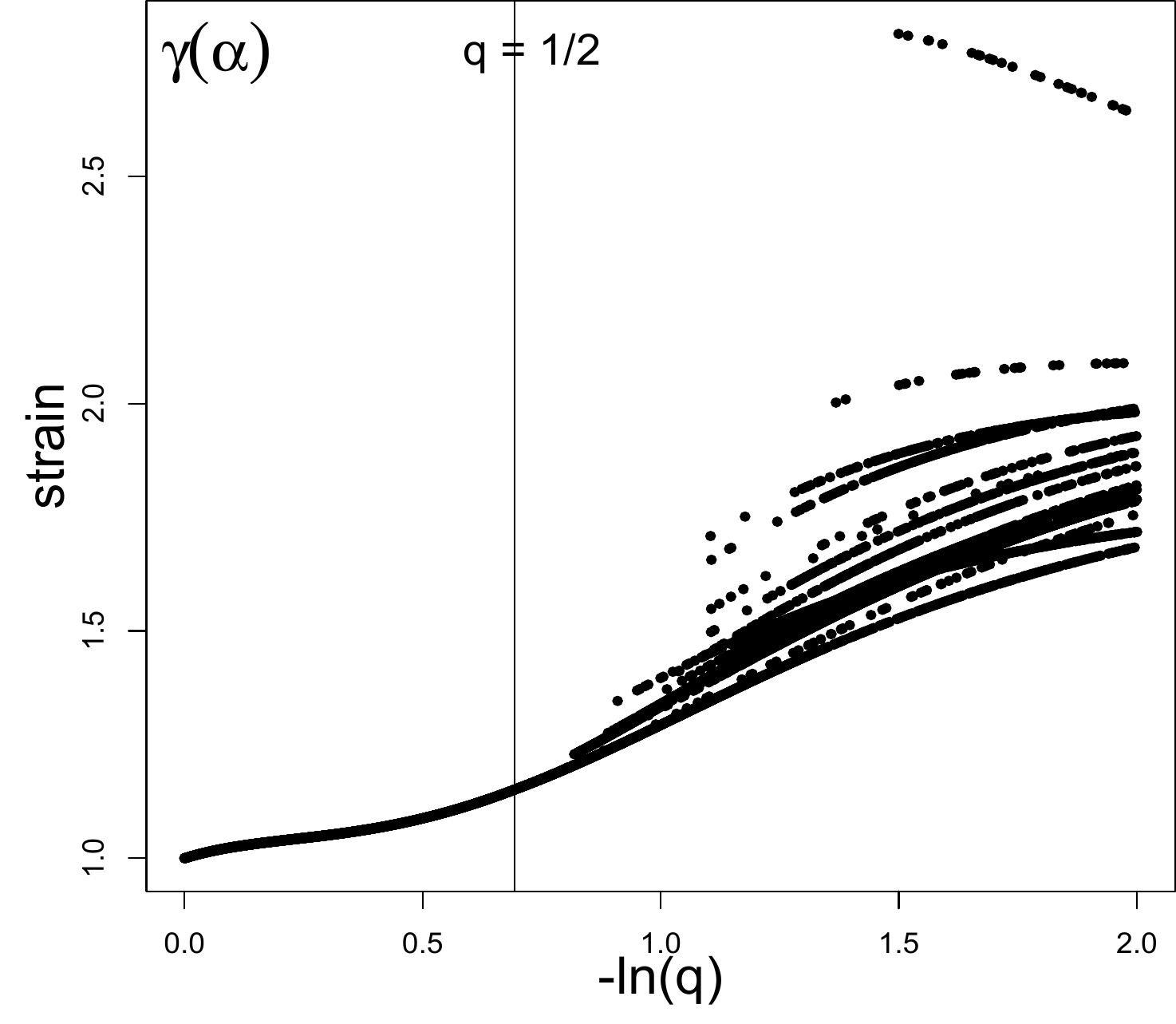}
\includegraphics[width=2.1in]{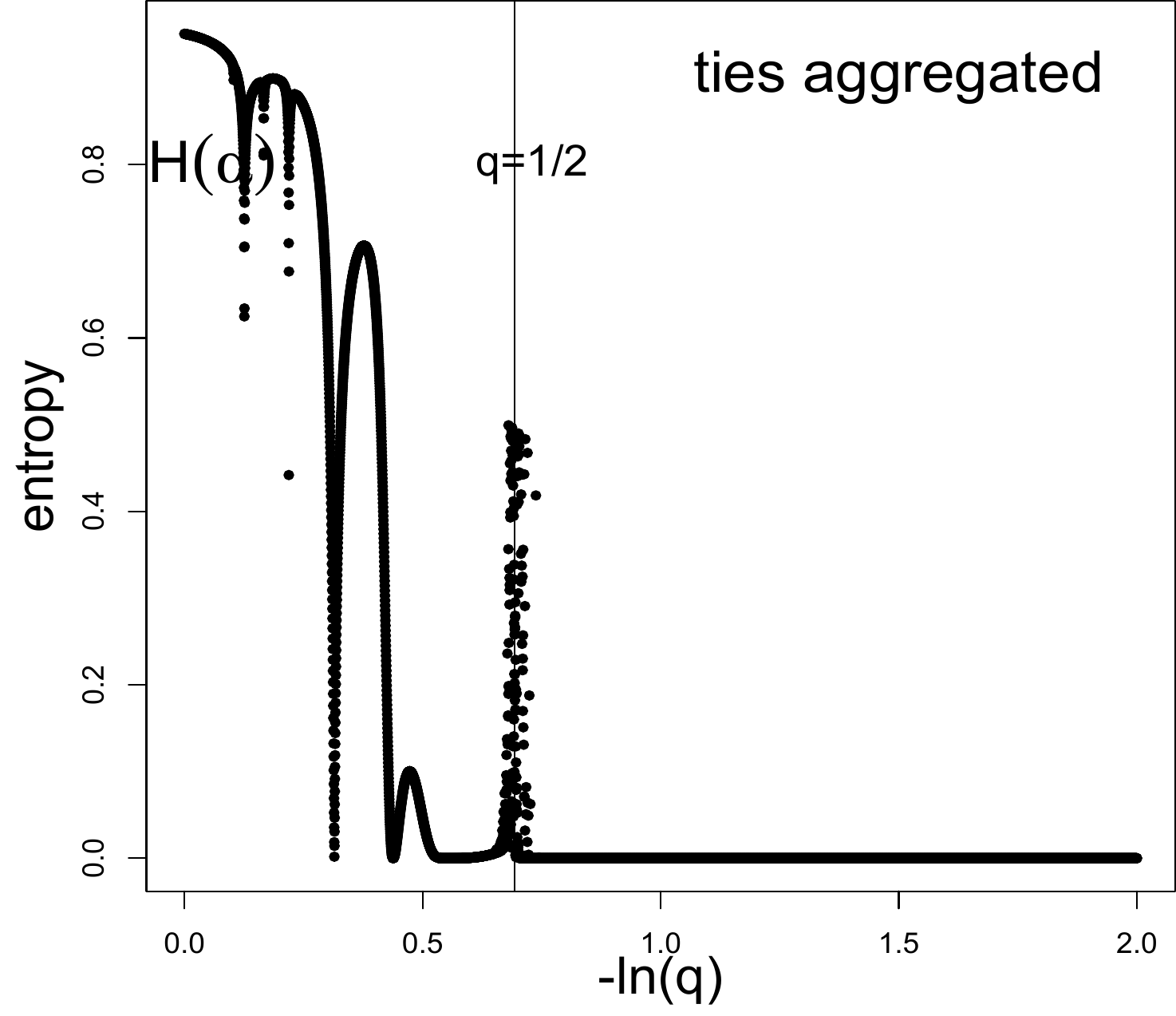}\hspace{0.5cm}
\includegraphics[width=2.1in]{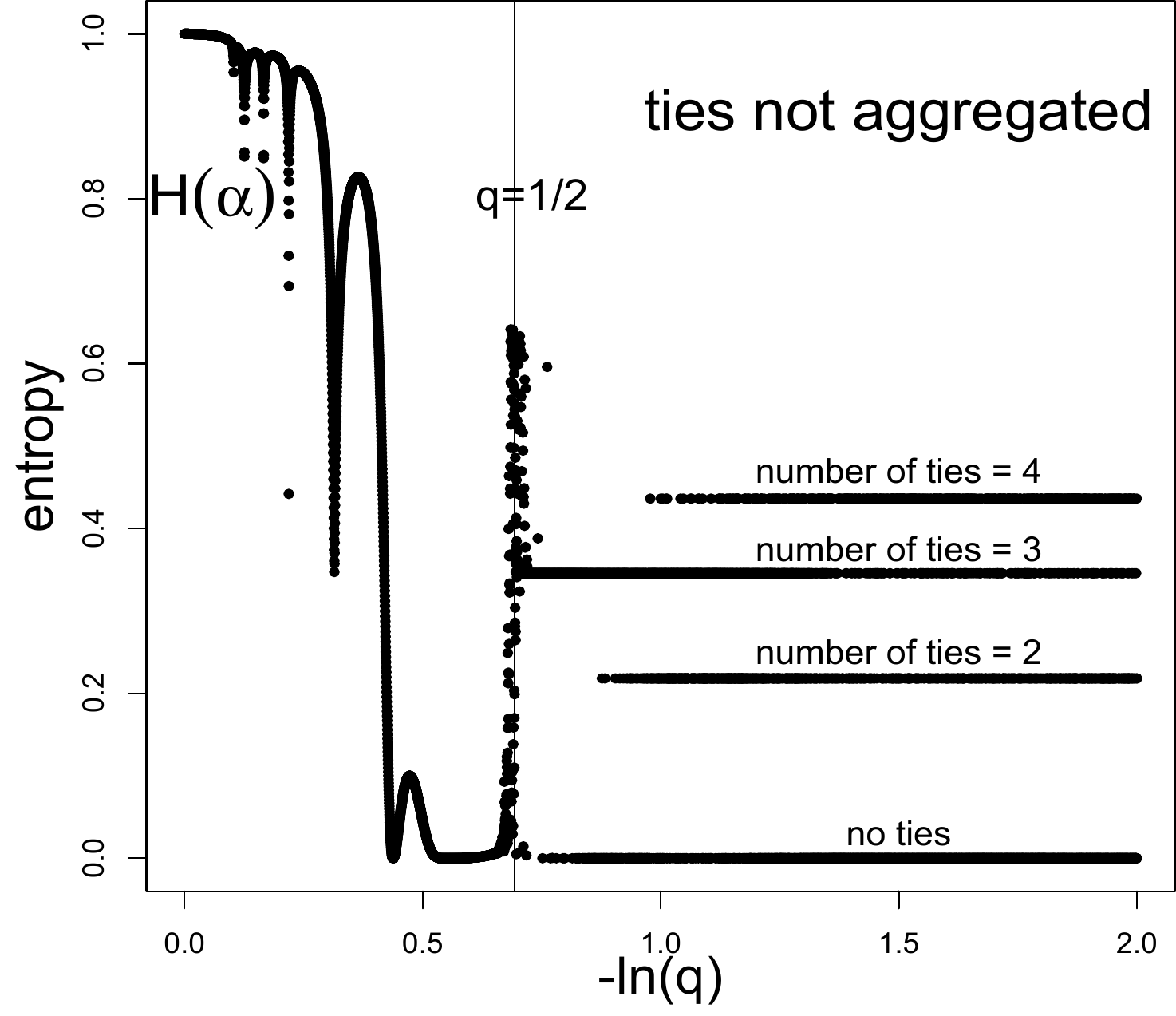}
\caption{Location estimates for the copper concentrations data, using the power transformation $\phi(D)=D^q$, for various values of $q$ (top left). Corresponding strain (top right) and entropies (bottom), exhibiting spurious, non-zero values when ties are not aggregated (bottom right).}
\label{figpower}
\end{center}
\end{figure}

Figures \ref{fourtransforms} and \ref{figpower} exhibit the values $a$ resulting from the iteration of equations (\ref{solmain}), for various transformations and various range of parameters, starting with an initial  value drawn as
 $a_0\sim U(2,6)$ (which de facto excludes the solutions around $x=28.95$, for readability sake). 
 
 \
 
The rectifiable transformations (Figure   \ref{fourtransforms}) yield distributed solutions (in accordance with 1a) and 1b) of Theorem \ref{tworegimes}), unique for large characteristic length (Theorem \ref{trecti}), and, in the case of bounded exponential and Tukey transformations, ``quasi-concentrated" for large characteristic length  (Theorem \ref{trsddsti}).

\

By contrast, the power transformation (Figure \ref{figpower}) exhibits a phase transition at $q=1/2$ (the median) between distributed and concentrated regimes (Section \ref{potrans} and Theorem  \ref{tworegimes}). The value of the strain can be shown to converge for $q\to0$ to $2(1-f_{i_0})/(1-\sum_j f_j^2)$, where $f_{i_0}$ is the weight of  the observation on which the minimum is concentrated; it ranges   from  $1.82$ (ties occurring four times) to $2.09$ (no ties).

As expected, the value of the entropy is zero in the concentrated regime. 
Figure \ref{figpower} bottom right demonstrates its deceptive behavior  when the 
 aggregation of ties is neglected. Also, the entropy dramatically decreases when the distributed solution crosses the sample values.

\section{The multivariate case}
\label{multi}
The robust estimation scheme (\ref{solmain}) depends on the $n\times n$ distance matrix only,  and not on the coordinates nor on their dimension. The multivariate nature of data is then of little  concern, in contrast to early attempts of generalizing $M$-estimates (such as the 
``peeling"  or the ``iterative trimming"    schemes  cited in Maronna (1976)). Rather, the success of the multidimensional implementation boils down to the ability to define an univocal squared Euclidean distance  for multidimensional data. This is typically the case for the chi-square metric of {\em Correspondence Analysis} (CA)  (Benz\'ecri 1973, Greenacre 1984),   illustrated below.

\subsection{Illustration: scientific collaborations of India}
 Consider the contingency table $N=(n_{ig})$ giving the number of scientific publications (1993-2000) co-authored by  India and   foreign countries $i=1,\ldots, 29$ (rows), sorted by  disciplines $g=1,\ldots,22$ (columns) (Anuradha and Urs 2007). The natural distance between countries is the {\em chi-square} distance 
\begin{equation}
\label{chis}
D_{ij}=\sum_g \rho_g (\frac{n_{ig}n_{\bullet \bullet}}{n_{i\bullet} n_{\bullet g}}-
\frac{n_{jg}n_{\bullet \bullet}}{n_{j\bullet} n_{\bullet g}})^2
\qquad\qquad
\rho_g=\frac{n_{ \bullet g}}{n_{\bullet \bullet}}   
\qquad\qquad
f_i=\frac{n_{i\bullet}}{n_{\bullet \bullet}}   
\end{equation}
where the dot symbol  sums over all values of the corresponding index. The associated inertia $\Delta(f)$ is, up to the total count,  the chi-square measure of rows-columns dependence; the determination of its  low-dimensional projection defines the well-known  Correspondence Analysis procedure.

\

\begin{figure}
\begin{center}
\includegraphics[width=2.3in]{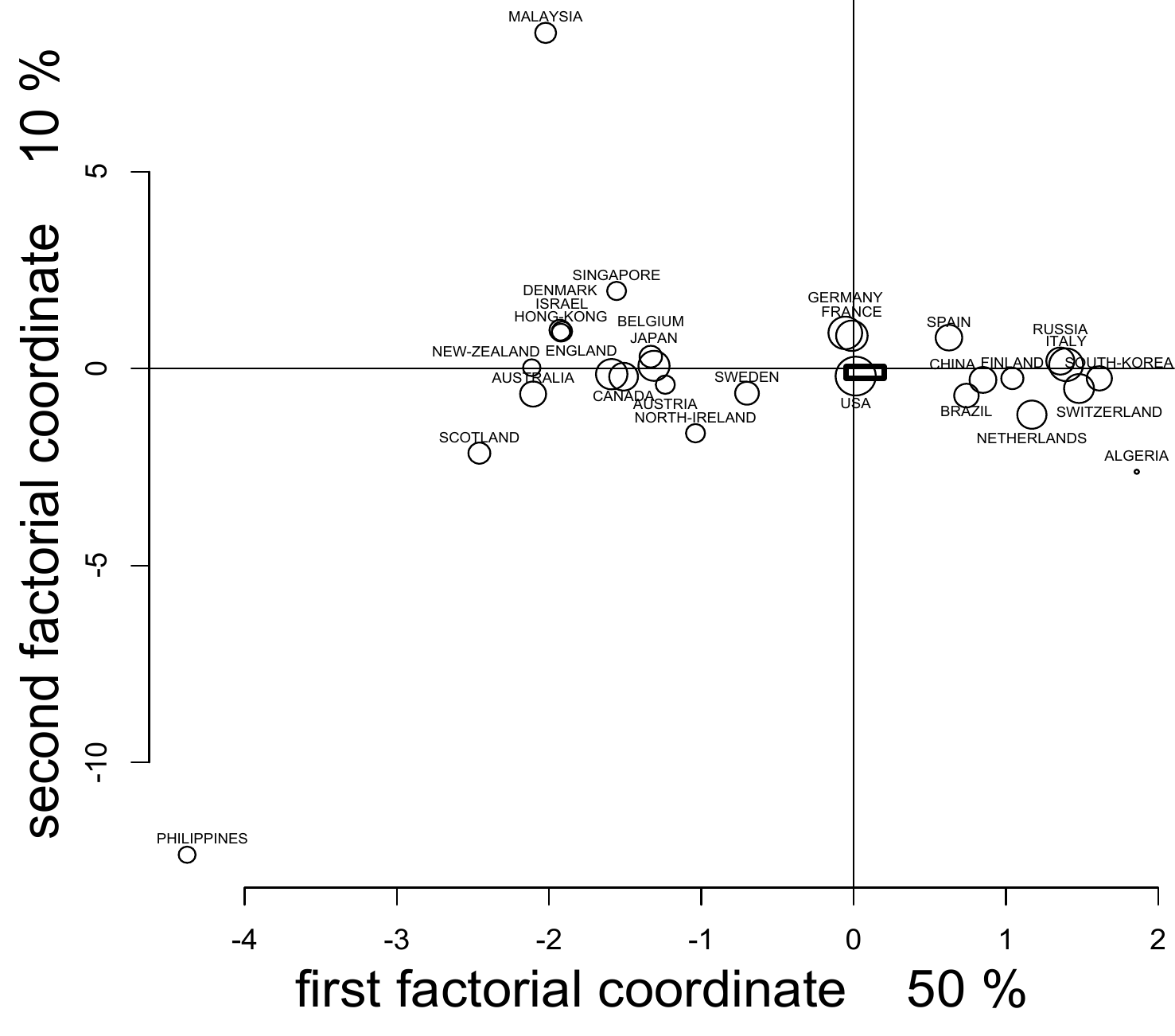}
\includegraphics[width=2.1in]{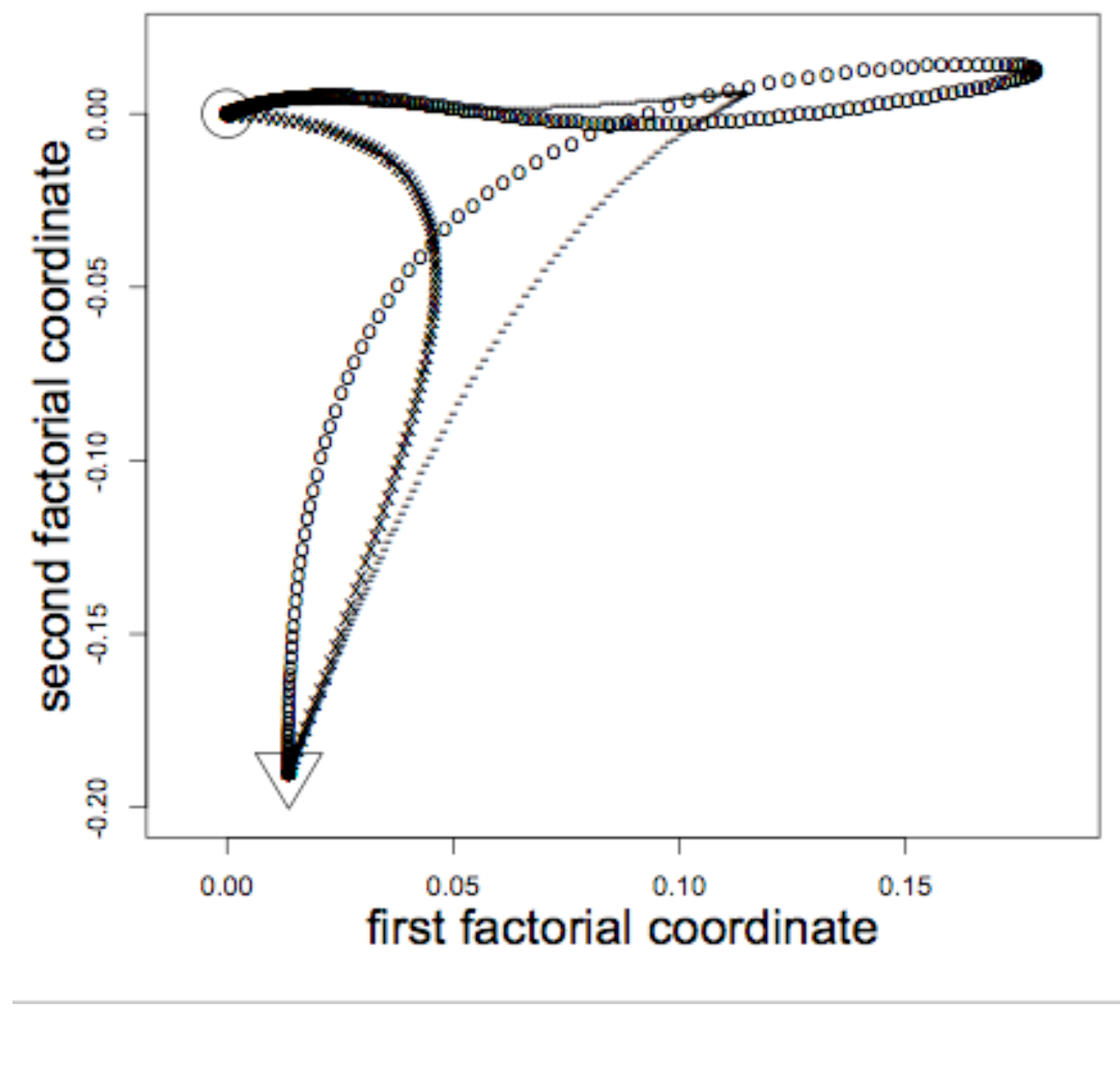}
\caption{Left: low-dimensional factorial coordinates of the countries, expressing 60\% of the inertia.  Right: zoom of the rectangle in the left figure, depicting the positions of the phi-estimate of location, under  power ({\tt{X}}), exponential ({\tt{O}}) and logarithmic  ({\tt{-}}) transformations, with parameters ranging in the distributed regime.}
\label{calacallo}
\end{center}
\end{figure}

 CA extracts uncorrelated factorial coordinates $x_{i\beta}$ such that 
 $D_{ij}=\sum_{\beta\ge0}(x_{i\beta}-x_{j\beta})^2$, each dimension $\beta$ accounting for a  proportion of  explained inertia, decreasing  with $\beta$.  Figure \ref{calacallo} left exhibits the countries coordinates in the first and second dimensions. The size of the circles depicts the countries weights, maximum for the
the USA, whose coordinates are also the closest to the origin.
 
 Figure \ref{calacallo} right depicts the plane trajectories, under various transformations in the distributed regime, of the centroid projections $\sum_i\alpha_i x_{i\beta}$ for $\beta=1,2$, where $\alpha$ is the profile resulting from (\ref{solmain}). 
 All  trajectories, diverse in shape and small in amplitude, lead from the
 the mean at the origin  $\bar{x}_{f\beta}=\sum_i  f_i x_{i\beta}=0$ (circle)
  to the coordinates of the USA (triangle).

\section{Conclusion}
Classical multidimensional scaling extracts the coordinates of the objects, uniformly weighted or not, from their Euclidean inter-distances only. In the same circumstances, this paper proposes a class of robust location estimates, based upon a Schoenberg transformation of distances, and  exhibiting a wide-ranging behavior. Further studies on their sensitivity, breakdown points and the like should presumably help  addressing the  question ``which transformation should be used in which context", left aside in the present set-up.

Additional guidance is also expected from a probabilistic point of view,  deliberately avoided as well in this paper. Manifestly, Schoneberg transformations can serve in  generating new densities: typically, replacing the distances in an univariate normal distribution by their square root produces an exponential distribution, which constitutes,  in the spirit of Section \ref{hdes}, a  normal distribution in some higher-dimensional embedded space. In that context, maximum likelihood considerations should also help in selecting a suitable family of transformations - whose Euclidean nature constitutes a favorable circumstance, regarding their tractabillity.

\section{Appendix}

\begin{proof}[of Theorem \ref{tworegimes}]
The derivative of a rectifiable transform is bounded. Setting to zero the derivative of  $\Gamma(a)$ with respect to the 
$l$-th component of the centroid $a_l=\sum_i\alpha_i x_{il}$ yields 
\begin{displaymath}
-2\sum_i f_i \varphi'(D_{ia})(x_{il}-a_l)=0
\end{displaymath}
which is the first identity in (\ref{solmain}), the second one resulting from the strong Huygens principle of section \ref{huyhuy}.  This solution is a minimum if, for all $t$
\begin{displaymath}
0\le\frac12\sum_{kl}\frac{\partial^2\Gamma(a)}{\partial a_k \partial a_l}t_kt_l
=\sum_i f_i [\varphi'(D_{ia}) \|t\|^2+2\varphi''(D_{ia})((x_{i}-a)'t)^2]
\end{displaymath}
Setting $t=b-a$ yields $((x_{i}-a)'t)^2=D_{ia}\|t\|^2\cos^2\theta_{ib}^a$, proving (\ref{scdo}). Equation  (\ref{stabily2bis}) follows form
$\phi''(D)\le0$ and $\cos^2\theta_{ib}^a\le 1$.
$\Box$\end{proof}

An alternative,  less direct proof also obtains by considering perturbing the profile elements $\alpha_i$  rather than the centroid components $a_l$. The second identity in (\ref{solmain}) yields
\begin{displaymath}
\frac{\partial D_{ia}}{\partial \alpha_j}=
D_{ij}-\sum_k\alpha_k D_{jk}=D_{ij}-D_{ja}-\Delta(\alpha)
\end{displaymath}
and hence the minimum condition
\begin{displaymath}
0=\frac{\partial\Gamma(a)}{\partial \alpha_j}=
\sum_i f_i \varphi'(D_{ia})(D_{ij}-D_{ja}-\Delta(\alpha))
\end{displaymath}
which, by comparison to the identity $\sum_i \alpha_i (D_{ij}-D_{ja}-\Delta(\alpha))=0$, yields the first identity in  (\ref{solmain}). 
Further derivating with respect to $\alpha_k$ yields together with (\ref{solmain})
the stability condition
\begin{eqnarray}    \sum_{jk}\frac{\partial^2 \Gamma(a)}{\partial \alpha_j\partial \alpha_k}z_jz_k& = & -\sum_i f_i \varphi'(D_{ia})\sum_{jk}D_{jk}z_jz_k \notag \\
& +&
\sum_i f_i \varphi''(D_{ia})[\sum_j (D_{ij}-D_{aj}-\Delta({\alpha}))z_j]^2\ge0
\label{truc}\end{eqnarray}
where $z$ is an admissible infinitesimal variation of the profile, that is of the form
 $z={\beta}-{\alpha}$ where ${\beta}$ is another profile.
But $\sum_{jk}D_{jk}z_jz_k=-2D_{ab}$ (e.g. Bavaud 2011) and 
\begin{displaymath}
\sum_j (D_{ij}-D_{aj}-\Delta({\alpha}))z_j= D_{ib}-D_{ab}-D_{ia}=
-2\sqrt{D_{ab}D_{ia}}\cos\theta_{ib}^a
\end{displaymath}
by the cosine formula.  $\Box$

\begin{proof}[of Theorem \ref{trecti}]
By hypothesis, $0<\varphi'(0)<\infty$. From   (\ref{schtr}),  
\begin{displaymath}
|D\phi''(D)|=\int_0^\infty \lambda D\exp(-\lambda D)g(\lambda)\: d\lambda\le \frac{1}{e}
\int g(\lambda)\: d\lambda=\frac{1}{e} \varphi'(0)
\end{displaymath}
implying  $\lim_{D\to0}D\phi''(D|\delta)=0$ by dominated convergence. Hence $\chi(0)>0$, and, by continuity  $\chi(D)>0$ for $D$  small enough, or equivalently for $\delta$ large enough. As a consequence, the stability condition 
(\ref{stabily2bis}) of Theorem (\ref{tworegimes}) holds irrespectively of the value of $a$, that is
  $\Gamma(a)$ is convex for  $\delta$ large enough, and  thus possesses an unique minimum. The latter   tends to $\bar{x}_f$ in the limit $\delta\to0$, in view of
$\varphi(D)=\frac{\eta(\delta)}{\delta} (D+0(D^2/\delta^2))$.  $\Box$
\end{proof}


 \end{document}